\documentclass[12pt]{article}

\usepackage[T2A]{fontenc}
\usepackage[utf8]{inputenc}

\usepackage{amsfonts}
\usepackage{amssymb}
\usepackage{amsmath}
\usepackage{amsthm}

\def \le {\leqslant}
\def \ge {\geqslant}

\theoremstyle{plain}

\begin{document}
\begin{Huge}
\centerline{\bf On the sequence $\alpha n!$}
\end{Huge}
\begin{Large}
\vskip+0.5cm
\centerline{Alena Aleksenko\thanks{Department of Mathematics, Aveiro University, Aveiro 3810, Portugal.  }}
\end{Large}
\vskip+2.0cm

\begin{small}
{\bf Abstract.} \,\, We prove that there exists $\alpha \in \mathbb{R}$ such that for any $N$  the dicrepancy $D_N$  of the sequence
$\{\alpha n !\} ,\, 1\le n \le N$ satisfies $ D_N = O(\log N)$.
\end{small}
\vskip+1.0cm

{\bf 1. Low discrepancy sequenses.}

Consider a sequence  $\xi_k, k =1,2,3...$ of points from the interval $[0,1)$.
 The discrepansy $D_N$ of the first $N$ points of the sequence
is defined as
$$
D_N = \sup_{\gamma\in [0,1)}  |N\gamma - \# \{ j :\,\,\, 1\le j \le N, \xi_j \le \gamma\}|.
$$
According to the famous W.M. Schmidt's theorem for any infinite sequence one has
$$
\liminf_{N\to \infty} \frac{D_N}{\log N} >0.
$$
This statement is sharp enough. Consider a real number $\alpha$ with bounded partial quotients in its continued fraction expansion. Then the discrepancy of the sequence
$$
\{ \alpha  k\} ,\,\,\ k =1,2,3,...
$$
(here $\{\cdot\}$ stands for the fractional part)
satisfied the inequality
$$
D_n \le M \log N,
$$
where $M$ depends on the bound  for the partial  quotients of $\alpha$. Roughly speaking, this result was obtained long ago by Ostrowski \cite{O} and Khintchine \cite{H}.
For further information concerning discrepancy bounds one cas see wonderful books \cite{T,KN} and \cite{MA}.

As for exponentionally  increasing sequences we would like to refer to Levin's paper \cite{L}. Given integer $q\ge 2$ Levin proved the existence of  real $\alpha$ such that the discrepancy $D_N$ for the sequence
\begin{equation}\label{exp}
\{ \alpha q^k\},\,\,\ k =1,2,3,..
\end{equation}
satisfies the inequality
\begin{equation}\label{exp1}
D_N  = O((\log N)^2).
\end{equation}
Up to now it is not known if there exists $\alpha$ such that the order of the discrepancy for the sequence (\ref{exp})
is smaller that  that from (\ref{exp1}). However if the sequence  $n_k$ grows faster than exponentially, it is quite easy to construct $\alpha$ such that the dicsrepancy of the sequence
$\{ \alpha n_k\}$ is small. This is just the purpose of the present short communication.

{\bf Theorem 1.}\,\, {\it
Suppose that a sequence of positive  numbers $ n_k, \, k=1,2,3,...$
 satisfy the condition
\begin{equation}\label{cond}
\inf_k \,   \frac{n_{k+1}}{k n_k}
>0.
\end{equation}
There exists $\alpha \in \mathbb{R}$ such that for any $N$  the dicrepancy $D_N$  of the sequence
$\{\alpha n_k\} ,\, 1\le k \le N$ satisfies $ D_N = O(\log N)$.}

{\bf Corolary 2.}\,\, {\it
 Then there exists $\alpha \in \mathbb{R}$ such that for any $N$  the dicrepancy $D_N$  of the sequence
$\{\alpha n !\} ,\, 1\le n \le N$ satisfies $ D_N = O(\log N)$.}

Korobov (\cite{K}, Theorem 3 and Example 1) constructed real numbers $\alpha$ for which the sequence $\{\alpha k!\}$ is uniformly distributed.
However his construction does not give optimal bounds for the discrepancy.

\vskip+0.3cm
{\bf 2.  Lemmas.}

In the sequel $F_i$ stands for the $i$-th Fibonacci number, so
$ F_0 = F_1 = 1, F_{i+1} =  F_i + F_{i-1}$ and
$ \phi = \frac{1+\sqrt{5}}{2}$.

{\bf Lemma 3.}\,\, {\it
Any positive integer $N$ can be represented in a form
$N = \sum_{i=1}^r b_i F_i$, where $ b_1 \in \{1,2,3\}, b_i\in \{ 1,2\}, 2\le i \le r$ and
$ r \le1+\log_\phi N$.}

Proof.  It is a well-known fact that  any positive integer
can be represented in a form  $N = \sum_{i=1}^t a_i F_i$ with $a_i \in\{0,1\}$,
$t \le 1+\log_\phi N$
and
in the sequence
\begin{equation}\label{seq}
a_1,a_1,...,a_{t}
\end{equation}
there is no two consequtive ones.
Now we give an algorithm how to construct from the sequence (\ref{seq}) a sequence $b_1,b_2,...,b_r$
with al positive
$b_i$ and $r\le t$.

We shall use two procedures.

{\bf Procedure 1.}   If we have two consecutive zeros, that is we have a pattern
$a_{i},0,0,a_{i+3}$, with $a_{i+3} =1$, we can replace it by the pattern
$a_{i},1,1,0$. The sum $N = \sum_{i=1}^t a_i F_i$ will remain the same as $ F_{i+3} = F_{i+2} + F_{i+1}$.
But the number of zeros in the sequence will decrease by one.

Procedure 1 enables one to get from the sequence (\ref{seq}) another sequence  of ones and zeros without  two consecutive zeros.

{\bf Procedure 2.}   If we have a pattern $a_i,0,a_{i+2}$  we may replace it by the pattern
$ a_{i}+1,1, a_{i+2}-1$. If $a_{i+2}= 1$ then there will be no zero in $(i+1)$-th position but a zero will appear in $(i+2)$-th  position. The total number of zeos wil not change.

The algorithm is as follows.
Procedure 1 enables one to get a sequence
\begin{equation}\label{seqs}
a_1',a_1',...,a_{t'}'
\end{equation}
 of ones and zeros where there in no two consecutive zeros
with $t' \le t$ and $ a_{t'} = 1$.

If $a_1'=0$ then $a_2' = 1$ and we may replace
the pattern $a_1',a_2'$ by $ 2,0$ as $ F_2 = 2= 2 F_1$.
So we may suppose that in (\ref{seqs}) one has $a_1' \in \{1,2\}$. But it may happen that $a_2'= a_3' =0$.
By applyinng procedures 1 we obtain a sequence
\begin{equation}\label{seqs}
a_1'',a_1'',...,a_{t''}''
\end{equation}
with $t''\le t'$ where
 $a_1' \in \{1,2\}$ and all other elements are equal to $ 1 $ or $0$ with no two consecutive zeros.
Now wew take the zero in the smallest position and apply Procedure 2. The  zero will turn into the next position. Then either there are two consecutive zeros  (and we can reduce the number of zeros by Procedure 1) or we  can move the zero in the next position again. Each "moving to the next position" increases the previous digit by 1.
But all the time the previous digit is 1. The only exception is in the very beginning of the process, when $a_1'' = 2$. Then  $a_1''$ must turn into 3.

In such a way we get the necessary representation for $N$. $\Box$.

From Lemma 3 we immediarely deduce

{\bf Corollary 4.}\,\,{\it
For any positive integer $ N$
the set
 of the first $N$ positive integers can be partitioned into segments of consecutive integers in the following way:
\begin{equation}\label{erepar}
\{ 1,2,...,N\} = {\cal A}
 \sqcup \,\bigsqcup_{i = 2}^r \,\bigsqcup_{j =1}^{b_i}\,
\{ R_{i,j}, R_{i,j}+1,...,R_{i,j}+F_i-1\},
\end{equation}
where ${\cal A} = \{1\}$ or $\{1,2\}$ or $\{1,2,3\}$,  $b_i \in \{ 1.2\}$  and
\begin{equation}\label{ere}
R_{i,j} \ge F_i,
\end{equation}
\begin{equation}\label{log}
r \le 1+\log_\phi N.
\end{equation}

}

Let $R, i $ be positive integers. We consider the sequence
\begin{equation}\label{fi}
\{
\phi k \},\,\,\  R\le k < R+F_i
\end{equation}
The following statement is well-known.
It is the main argument of the classical proofs
of the
logarithmic order of discrepansy  of the sequence $\{\alpha k\} $, in the case when $\alpha$ has bounded partial quotinents in its continued fraction expansion.
It immediately follows from the inequality
$||\phi F_i||\le 1/F_i$, where  $||\cdot ||$ stands for the distance to the nearest integer.
 It means that the the set (\ref{fi}) is close to the set
$$
\frac{\nu}{F_i},\,\,\, 1,\,\,\,0\le \nu \le F_i -1,
$$
and hence discrepancy of the sequence (\ref{fi}) is bounded.

{\bf Lemma 5.}\,\, {\it There is a substitution $ \sigma_1,...,\sigma_{F_i}$ of the sequence $ 1,...,F_i$ such that
$$
\left|
\{
\phi (R+ k) \} -  \frac{\sigma_k}{F_i}
\right| \le \frac{1}{F_i},\,\,\, 0\le k \le F_i-1.
$$
}

{\bf Lemma 6.}\,\, {\it Consider an arbitrary sequence $\xi_k , \, k=1,2,3,...$ from the interval $[0,1)$.
Suppose that a sequence $n_k$ satisfy (\ref{cond}). Then there exist $c>0$ and  $\alpha\in \mathbb{R}$ such that
 \begin{equation}\label{vovo}
||\alpha n_k - \xi_k ||\le \frac{c}{k},\,\,\ k =1,2,3,... .
\end{equation}
 .
}

Proof. From (\ref{cond}) we see that
$ \frac{n_{k+1}}{n_k}\ge {\kappa}{k}$ for  some positive $\kappa$.
Fix $k$. Then the set
$$
\left\{ \alpha \in \mathbb{R}:\,\,\,  || \alpha n_k -\xi_k||\le \frac{c}{k}\right\}
$$
is a union of segments of the form $\left[ \frac{\xi_k + z}{n_k} - \frac{c}{kn_k}, \frac{\xi_k + z}{n_k} +\frac{c}{kn_k}\right], z \in \mathbb{Z}$.  The length of each segment is equal to $\frac{2c}{kn_k}$.
The distance between the centers of neighbouring segments is equal to $ \frac{1}{n_k}$.
We see that $ \frac{c}{kn_k} \ge \frac{c\kappa}{n_{k+1}}$.
So if $c$ is large enough, we  can choose integers $z_k$ to get a seuence of nested segments
$$
\left[ \frac{\xi_1 + z_1}{n_1} - \frac{c}{n_1}, \frac{\xi_1 + z_1}{n_1} +\frac{c}{n_1}\right] \supset ... \supset
\left[ \frac{\xi_k + z_k}{n_k} - \frac{c}{kn_k}, \frac{\xi_k + z_k}{n_k} +\frac{c}{kn_k}\right]
\supset
$$
$$
\supset
\left[ \frac{\xi_{k+1} + z_{k+1}}{n_{k+1}} - \frac{c}{(k+1)n_{k+1}}, \frac{\xi_{k+1} + z_{k+1}}{n_{k+1}} +\frac{c}{(k+1)n_{k+1}}\right]
\supset ... .
$$
The common point of these segments satisfies (\ref{vovo}).$\Box$

\vskip+0.3cm
{\bf 3. Proof of Theorem 1.}

We take the sequence
$\xi_k =\{ \phi k\}, k =1,2,3,...$ and apply Lemma  6. Then we get real $\alpha$. This is just the number what we need.
Take positive integer $N$. Then in the decomposition (\ref{erepar}) for each segment
$\{ R_{i,j}, R_{i,j}+1,...,R_{i,j}+F_i-1\}$ its right endpoint is $\ge F_i$. So  from the inequalities of Lemmas 5,6 we get
$$
\left|\left|
\alpha n_k -  \frac{\sigma_k}{F_i}\right|\right|\le
\frac{1}{F_i}+\frac{c}{R_{i,j}}\le \frac{1+c}{F_i},\,\,\,\forall i.j.
$$
This means that each sequence
\begin{equation}\label{por}
\{\alpha n_k\},\,\,\,  R_{i,j}< k \le R_{i,j} +F_i
\end {equation}
has discrepancy $O(1)$.
But the sequence  $\{\alpha n_k\}, \, 1\le k \le N$ is partitioned into $ O(\log N)$ sequences of the form (\ref{por}), as all $b_i $ are bounded by  3 and we have estimate (\ref{log}).
So this sequence has discrepancy $O(\log N)$.$\Box$

\textbf{Acknowledgment.} Author is very grateful to N.Moshchevitin for his permanent help and support. This work was supported by {\it FEDER} founds through {\it COMPETE}--Operational Programme Factors of Competitiveness (``Programa Operacional Factores de Competitividade'') and by Portuguese founds through the {\it Center for Research and Development in Mathematics and Applications} (University of Aveiro) and the Portuguese Foundation for Science and Technology (``FCT--Fund\c{c}\~{a}o para a Ci\^{e}ncia e a Tecnologia''), within project PEst-C/MAT/UI4106/2011 with COMPETE number FCOMP-01-0124-FEDER-022690, and by the FCT research project PTDC/MAT/113470/2009


\begin{thebibliography}{100}


\bibitem{T}
 M. Drmota, R.  Tichy,\,\,
Sequences, discrepancies, and applications, ?ect. Not. Math., 1651, 1997.


\bibitem{H}
A. Khintchine,\,\,
Ein Sats \"uber Kettenbr\'uche, mit  arithmetischen Anwendungen, Math. Zeitschrift 18, no. 3,4  (1923), 289 - 306.

\bibitem{K}
N.M. Korobov,\,\,
\"Uber einige Fragen der Gleichverteilung. (in Russian)
 Izv. Akad. Nauk SSSR, Ser. Mat. 14, 215-238 (1950).

\bibitem{KN}
L. Kuipers, Y. Niederreiter,\,\,
Uniform Distribution of Sequences, Dovers Publication. 1974.

\bibitem{L}
M. B. Levin,\,\,
On the discrepancy estimate of normal numbers.
 Acta Arith. 88, No.2  (1999), 99 - 111.

\bibitem{MA}
J. Matousek,\,\, Geometric Discrepancy,
Algorithms and Combinatorics
Volume 18,
Springer, 2010.


\bibitem{O}
A. Ostrowski,\,\,
Bemerkungen zur Theorie der diophantischen Approximationen.
Hamb. Abh. 1 (1921), 77 - 98.

\bibitem{S}
W.M. Schmidt,\,\,
On irregularities of distribution VII. Acta
Arith., 21 (1972), 45–50.


\end{thebibliography}
\end{document}